# A Note on Pépin's counter examples to the Hasse principle for curves of genus 1

Franz Lemmermeyer

March 18, 1999

**Abstract**

In a series of articles published in the C.R. Paris more than a century ago, T. Pépin announced a list of "theorems" concerning the solvability of diophantine equations of the type $ax^4 + by^4 = z^2$. In this article, we show how to prove these claims using the structure of 2-class groups of imaginary quadratic number fields. We will also look at a few related results from a modern point of view.

## 1 Introduction

It is well known that curves of genus 0 defined over $\mathbb{Q}$ satisfy the Hasse principle: they have a rational point if and only if they have a $\mathbb{Q}_p$-rational point for every completion $\mathbb{Q}_p$ of $\mathbb{Q}$. It is similarly well known that the Hasse principle fails to hold for curves of genus 1, the first counter example $2z^2 = x^4 - 17y^4$ being due to Lind [11] and Reichardt [22].

In a series of articles [15, 17, 18, 19], Théophile Pépin announced 93 theorems asserting that certain equations of the type $aX^4 + bY^4 = Z^2$ were not solvable in integers (nontrivially, that is). In order to get nontrivial results, Pépin looks at equations whose underlying conics $ax^2 + by^2 = z^2$ do have rational solutions (see [15]):

> Les cas où l'équation indéterminée $aX^4 + bY^4 = Z^2$ n'admet pas de solution rationnelle sont fort nombreux, même quand l'équation $ax^2 + by^2 = z^2$ est résoluble en nombres entiers. Néanmoins on ne connait encore qu'un petit nombre de théorèmes sur ce sujet.

He then starts listing his results without proof, and among his examples there are some that claim the nonexistence of rational points on some curves of genus 1 that are everywhere locally solvable (Pépin probably wasn't familiar with the genus of algebraic curves, $p$-adic integers were invented around 1897 by Hensel, andthe first versions of the Hasse principle were discovered by Hasse in 1920, so the title of this paper is somewhat misleading). To the best of my knowledge no proofs for Pépin's claims have been supplied yet. The proofs that we give in Section 2 below are based on connections with the 2-class groups of complex quadratic number fields; such a connection is not unexpected in view of earlier work of Rédei [20] on similar diophantine problems.

## 2  A Few of Pépin's Results

Let us start with the following assertion taken from [17]:

**Proposition 1.** *Let $p$ be a prime of the form $p = 5m^2 + 4mn + 9n^2$; then the equation $px^4 - 41y^4 = z^2$ does not have rational solutions.*

*Proof.* First observe that we may assume that $x, y, z$ are integers; moreover, if $q$ is a prime dividing $x$ and $y$ then $q^2 \mid z$ since $41p$ is squarefree, hence we may assume that $x$ and $y$ are coprime. Since any common prime divisor of $x$ and $z$ divides $y$, we also see that $(x, z) = 1$, and similarly $(y, z) = 1$.

Now write the equation in the form $px^4 = N(z + y^2\sqrt{-41})$, where $N$ denotes the norm from the quadratic field $k = \mathbb{Q}(\sqrt{-41})$. It is easily seen that $x$ is always odd and that either $y$ or $z$ is even. In particular, the ideals $(z + y^2\sqrt{-41})$ and $(z - y^2\sqrt{-41})$ are coprime. This implies that $(z + y^2\sqrt{-41}) = \mathfrak{p}\mathfrak{a}^4$, where $\mathfrak{p}$ denotes a prime ideal above $p$ in $k$: in particular, the ideal class of $\mathfrak{p}$ is a 4th power.

On the other hand we have $p = 5m^2 + 4mn + 9n^2$. This implies $5p = (5m + 2n)^2 + 41n^2$, hence $\mathfrak{p}$ is in the same ideal class as one of the primes above 5. Now a simple computation shows that each prime $\mathfrak{q}$ above 5 generates an ideal class of order 4 (note that $5^4 = 16^2 + 3^2 \cdot 41$); since $\mathrm{Cl}(k)$ is cyclic of order 8, the class $[\mathfrak{q}]$ is not a fourth power: contradiction.   □

Note that the proof shows more than claimed since we may replace $y^4$ by $y^2$; the same remark applies to other claims of Pépin. Note that not all of Pépin's claims produce counter examples to the Hasse principle for curves of genus 1: for example, although the conic $px^2 - 41y^2 = z^2$ has a rational point whenever $p$ is as in Proposition 1, the equation $px^4 - 41y^4 = z^2$ is solvable 2-adically if and only if $p \equiv 1 \bmod 8$; since solvability at the other completions is easy to verify, this curve is a counter example to the Hasse principle if and only if $p \equiv 1 \bmod 8$. Here is an example from [15] whose proof uses ray class groups:

**Proposition 2.** *Let $p$ be a prime of the form $p = 9a^2 + 4b^2$; then the equation $px^4 - 36y^4 = z^2$ does not have rational solutions.*

*Proof.* As above, $px^4 = N(z + 6y^2 i)$ implies that the class of the prime ideal $\mathfrak{p}$ above $p$ in $k = \mathbb{Q}(i)$ is a fourth power in the ray class group $\mathrm{Cl}_k\{6\}$ modulo 6 of $k$. On the other hand, $p = 9a^2 + 4b^2$ implies that $[\mathfrak{p}]$ has order 2 in $\mathrm{Cl}_k\{6\} \simeq \mathbb{Z}/4\mathbb{Z}$, and this is a contradiction.   □

Below, we will give two more proofs of Proposition 2 (not because the result is so important but in order to illustrate the different techniques): in Section 3, we employ the arithmetic of quadratic number fields, whereas the proof in Section 4 uses nothing beyond unique factorization in integers.

In several reviews (see e.g. [14]) it was noticed that Pépin's examples all had the following form: $p$ is a prime of the form $\alpha a^2 + 2\beta ab + \gamma b^2$, $m = \alpha\gamma - \beta^2$ is a positive integer, and $px^4 - my^4 = z^2$ is the equation that is claimed to have no nontrivial solutions. In this generality, however, the theorem does not hold; moreover it includes examples of equations such that even the underlying conic does not have rational points: take for example $\alpha = 2$, $\beta = 0$ and $\gamma = 3$: then $m = 6$, and in fact the equation $px^4 - 6y^4 = z^2$ has no solutions if $p = 2a^2 + 3b^2$. It is, however, easy to verify that $px^2 - 6y^2 = z^2$ is not solvable 2-adically.

It is possible to get families of examples with solvable conics; the following contains several of Pépin's claims:

**Theorem 3.** *Assume that $p = \alpha^2 a^2 + 2\beta ab + \gamma b^2$ is a prime, and put $m = \alpha^2\gamma - \beta^2$. Then the conic $px^2 - my^2 = z^2$ has the rational point $(x, y, z) = (\alpha, b, \alpha a + \beta b)$, hence infinitely many.*

*If, in addition, $m \equiv 1 \bmod 8$ and $\alpha \equiv 3 \bmod 4$ are prime, then the equation $px^4 - my^4 = z^2$ does not have nontrivial rational solutions.*

*Proof.* If a conic defined over $\mathbb{Q}$ has a rational point $P$, then any line through $P$ with a rational slope $t$ will intersect the conic in another rational point, thus producing a rational parametrization of the conic.

Assume that $px^4 - my^4 = z^2$ is solvable, and that $x$, $y$ and $z$ are pairwise coprime integers. Since $(-m/p) = (-m/\alpha) = 1$, we find that $p\mathcal{O}_k = \mathfrak{p}\mathfrak{p}'$ and $\alpha\mathcal{O}_k = \mathfrak{a}\mathfrak{a}'$ split in $k = \mathbb{Q}(\sqrt{-m})$. Now $px^4 = N(z + y^2\sqrt{-m})$, and since $x$ must be odd we deduce that $(z + y^2\sqrt{-m}) = \mathfrak{p}\mathfrak{b}^4$ for some ideal $\mathfrak{b}$. This shows that the ideal class of $\mathfrak{p}$ is a fourth power.

On the other hand, $p = \alpha^2 a^2 + 2\beta ab + \gamma b^2$ gives $p\alpha^2 = (\alpha^2 a + \beta b)^2 + mb^2$, hence $\mathfrak{p} \sim \mathfrak{a}^2$. Since $m \equiv 1 \bmod 4$, one of the genus characters of $k/\mathbb{Q}$ is the nontrivial character modulo 4. Now genus theory implies that the ideal classes of prime ideals above primes $\equiv 3 \bmod 4$ are not squares: since $\alpha$ is such a prime, $[\mathfrak{a}]$ is not a square, and we have a contradiction since the 2-class group of $k$ is cyclic of order divisible by 4 (this consequence of the fact that $m \equiv 1 \bmod 8$ is prime is classical: see e.g. Rédei & Reichardt [21]). $\square$

Let us also mention another example of Pépin's: in [18] he claimed that $px^4 - 32y^4 = z^2$ does not have nontrivial solutions if $p = 4u^2 + 4uv + 9v^2$ (our proof uses the ray class group modulo 4 of $\mathbb{Q}(\sqrt{-2})$ and is left to the reader); writing this as $p = (2u + v)^2 + 8v^2$, we see that solvability of $px^4 - 2y^4 = z^2$ for primes $p \equiv 1 \bmod 8$ implies that $p = A^2 + 32B^2$. See Rose [24] for more on the corresponding elliptic curve $y^2 = x(x^2 - 2p)$.

Finally we would like to show that one of the curves studied by Lind and Reichardt can be treated with our method: since fourth powers of odd integers $\equiv 1 \bmod 16$, the equation $-2z^2 = x^4 - 17y^4$ implies that $4 \mid z$. But then $17y^4 = x^4 + 32z^2$, and since the ray class group modulo 4 of $\mathbb{Q}(\sqrt{-2})$ is cyclic of order 4, the rest is clear.

## 3   Modern Interpretation

Today we recognize an equation of type $pX^4 - mY^4 = Z^2$ as a torsor of an elliptic curve. In fact, consider any elliptic curve $E : y^2 = x(x^2 + ax + b)$ with a rational point $(0,0)$ of order 2. Then there exists another curve $\widehat{E} : y^2 = x(x^2 - 2ax + a^2 - 4b)$ and 2-isogenies $\phi : E \longrightarrow \widehat{E}$ and $\psi : \widehat{E} \longrightarrow E$ such that their composition is multiplication by 2. The torsors one needs to study for computing the Selmer group $S^{(\psi)}(\widehat{E}/\mathbb{Q})$ attached to $\phi$ are given by

$$\mathcal{T}^{(\psi)}(b_1) : N^2 = b_1 M^4 + a M^2 e^2 + b_2 e^4,$$

where $b_1$ is squarefree and $b_1 b_2 = b$. Define a group structure on such torsors by identifying them with the classes $b_1 \mathbb{Q}^{\times 2} \in \mathbb{Q}^\times / \mathbb{Q}^{\times 2}$; then the Selmer group $S^{(\psi)}(\widehat{E}/\mathbb{Q})$ consists of all torsors $\mathcal{T}^{(\psi)}(b_1)$ that are everywhere locally solvable. Its subgroup of torsors that have a rational point will be denoted by $W(\widehat{E}/\mathbb{Q})$, and the quotient group $\text{III}(\widehat{E}/\mathbb{Q})[\psi] := S^{(\psi)}(\widehat{E}/\mathbb{Q})/W(\widehat{E}/\mathbb{Q})$ is called the Tate-Shafarevich group of $\widehat{E}(\mathbb{Q})$ attached to $\psi$; note that this group is nontrivial if and only if one of the torsors $\mathcal{T}^{(\psi)}(b_1)$ is a counter example to the Hasse principle. The corresponding groups attached to $\phi$ are defined similarly. Computing the Selmer groups attached to $\phi$ and $\psi$ is called performing a first simple 2-descent; ultimatively, of course, one is interested in $W(E/\mathbb{Q})$ and $W(\widehat{E}/\mathbb{Q})$: a formula due to Tate states that the product of their cardinalities equals $2^{r+2}$, where $r$ is the $\mathbb{Z}$-rank of $E(\mathbb{Q})$ and $\widehat{E}(\mathbb{Q})$.

In this section we will study the elliptic curve $E : y^2 = x(x^2 - 4pq^2)$ using such a first 2-descent, where $p \equiv 1 \bmod 8$ and $q \equiv 3 \bmod 4$ are primes such that $(p/q) = +1$ (for $q = 3$, this is the curve occurring in Proposition 2). We will use the notation from [10]. The curve $E$ is 2-isogenous to $\widehat{E} : y^2 = x(x^2 + pq^2)$, and it is easy to see that $S^{(\phi)}(E/\mathbb{Q}) = \langle p\mathbb{Q}^{\times 2}\rangle = W(E/\mathbb{Q})$. A simple calculation reveals that $S^{(\psi)}(\widehat{E}/\mathbb{Q}) = \langle -\mathbb{Q}^{\times 2}, 2\mathbb{Q}^{\times 2}, p\mathbb{Q}^{\times 2}, q\mathbb{Q}^{\times 2}\rangle$, and that $W(\widehat{E}/\mathbb{Q}) \supseteq \langle -p\mathbb{Q}^{\times 2}\rangle$.

Next we consider some torsors in detail; since $-p\mathbb{Q}^{\times 2} \in W(\widehat{E}/\mathbb{Q})$, it is sufficient to look at $\mathcal{T}^{(\psi)}(b_1)$ for $b_1 \in \{p, \pm 2, \pm q, \pm 2q\}$.

The following (partial) result can be proved by elementary number theory alone:

**Proposition 4.** *If the torsor $\mathcal{T}^{(\psi)}(b_1)$ has a nontrivial rational point, then the conditions in the following table are satisfied:*

| $b_1$ | $-pb_1$ | condition |
|------:|--------:|-----------|
| $2$   | $-2p$   | $(2/p)_4 = +1$ |
| $-2$  | $2p$    | $(2/p)_4 = +1$ |
| $q$   | $-pq$   | $(q/p)_4 = +1$ |
| $-q$  | $pq$    | $(q/p)_4 = +1$ |
| $2q$  | $-2pq$  | $(2q/p)_4 = +1$ |
| $-2q$ | $2pq$   | $(2q/p)_4 = +1$ |

The proof proceeds case by case:

- $b_1 = q$ Here we find $qN^2 = M^4 - 4pe^4$; this implies $2 \mid M$, hence $qn^2 = 4m^4 - pe^4$, which gives $(q/p)_4 = (2n/p) = (n/p) = (p/n) = +1$.

- $b_1 = -q$ Now $-qN^2 = M^4 - 4pe^4$; here $M$ is odd, hence $(-q/p)_4 = (N/p) = (p/N) = +1$. Note that $(-1/p)_4 = 1$ since $p \equiv 1 \bmod 8$.

- $b_1 = 2pq$ Here we get $2qn^2 = pM^4 - e^4$; hence $(2q/p)_4 = (n/p) = (n'/p) = (p/n') = 1$.

- $b_1 = -2pq$ Here $-2qn^2 = pM^4 - e^4$, and as above we find $(2q/p)_4 = 1$.

- $b_1 = 2$ Here we start with $N^2 = 2M^4 - 2pq^2e^4$. Assume first that $q \equiv 3 \bmod 4$, that is, $(2/q) = -1$. Reducing modulo $q$ shows that we must have $q \mid M$. Put $N = 2qn$ and $M = qn$; then $2n^2 = q^2M^4 - pe^4$ gives $(2/p)_4 = (qn/p)$. But $(q/p) = 1$ and $(n/p) = (n'/p) = (p/n') = 1$, hence $(2/p)_4 = 1$. If $q \equiv 7 \bmod 8$, on the other hand, then $q \nmid M$, and we get $2n^2 = M^4 - pq^2e^4$, but the conclusion $(2/p)_4 = 1$ stays the same.

- $b_1 = -2$ Then $N^2 = -2M^4 + 2pq^2e^4$. In this case, $q \mid M$ if $q \equiv 7 \bmod 8$, and $q \nmid M$ otherwise. Assume $q \equiv 7 \bmod 8$ first. Then $N = 2qn$, $M = qm$ and thus $-2n^2 = q^2m^4 - pe^4$. This implies $(2/p)_4 = 1$ as usual, and the other case is treated similarly.

A similar consideration of $\mathcal{T}^{(\psi)}(p)$ produces no result, although we know by Pépin's result in the special case $q = 3$, plus the fact that $p = 9a^2 + 4b^2$ is equivalent to $(-3/p)_4 = -1$, that solvability implies $(3/p)_4 = 1$. So we better have a second look at our torsor $\mathcal{T}^{(\psi)}(p)$.

Here $N^2 = pM^4 - 4q^2e^4$, where $MN$ is odd and $e$ is even; we also know that $q \nmid M$ since otherwise reduction modulo $q$ would imply that $(-1/q) = +1$. Then $pM^4 = N^2 + 4q^2e^4 = (N + 2qe^2i)(N - 2qe^2i)$ and thus $N + 2qe^2i = \pi\mu^4$, where $\pi \in \mathbb{Z}[i]$ is a prime $\equiv 1 \bmod 2$ and $\mu\overline{\mu} = M$. Subtracting its conjugate from this equation gives $4qe^2i = \pi\mu^4 - \overline{\pi}\,\overline{\mu}^4$, and reducing modulo $\overline{\pi}$ shows that $(q/p)_4(-1/p)_8(e/p) = [\pi/\overline{\pi}]_4$. Now $[\pi/\overline{\pi}]_4 = (-4/p)_8$ (see Lemma 7 below) and $(e/p) = (e'/p) = (p/e') = 1$ (where $e = 2^j e'$ for some odd $e'$) implies that $(q/p)_4(-1/p)_8 = (-4/p)_8$, hence $(2q/p)_4 = +1$.

Now we have found a necessary condition but it isn't the one we were expecting. So let's have another try and factor the torsor over $k = \mathbb{Q}(\sqrt{p})$:

$$-q^2e^4 = \Big(\frac{N + M^2\sqrt{p}}{2}\Big)\Big(\frac{N - M^2\sqrt{p}}{2}\Big).$$

Assume for the moment that $k$ has class number 1. Then we get $N + M^2\sqrt{p} = 2\varepsilon\lambda^2\alpha^4$, where $\varepsilon$ is some unit in $\mathcal{O}_k$ and $N\lambda = q$. Taking the norm of both sides shows that $N\varepsilon = -1$, so up to squares (which we may subsume into $\lambda$) we have $\varepsilon = \pm\varepsilon_p$, where $\varepsilon_p > 1$ is the fundamental unit of $k$. We see that $N + M^2\sqrt{p} > 0$ under the embedding that takes $\sqrt{p}$ to the positive real square root of $p$, soe we must have $N + M^2\sqrt{p} = 2\varepsilon_p\lambda^2\alpha^4$. Subtracting this equation from its conjugate yields $M^2\sqrt{p} = \varepsilon_p\lambda^2\alpha^4 + \overline{\varepsilon}_p\overline{\lambda}^2\overline{\alpha}^4$, and now reduction modulo $\overline{\lambda}$ gives $[\varepsilon_p \sqrt{p}/\overline{\lambda}] = 1$. But Kummer theory and a few arguments about ramification (see e.g. [9]) show that $\mathbb{Q}(\sqrt{\varepsilon_p \sqrt{p}}\,)$ is the quartic subfield of the cyclotomic field $\mathbb{Q}(\zeta_p)$, hence $[\varepsilon_p \sqrt{p}/\overline{\lambda}] = (q/p)_4$, which is exactly what we wanted.

For the general case, we need a lemma:

**Lemma 5.** *Let $p$ and $q$ be odd primes and $h \geq 1$ an odd integer such that $r^2 - ps^2 = q^h$ with $2r, 2s \in \mathbb{Z} \setminus q\mathbb{Z}$. Then $(r + s\sqrt{p}\,)^h - (r - s\sqrt{p}\,)^h = 2S\sqrt{p}$, and $(S/q) = (s/q)$.*

*Proof.* We have $S = \binom{h}{1}r^{h-1}s + \binom{h}{3}r^{h-3}s^3 p + \ldots + \binom{h}{h}s^h p(h-1)/2$. Since $r^2 \equiv ps^2 \bmod q$, this implies that $S \equiv s^h p^{(h-1)/2}[\binom{h}{1} + \binom{h}{3} + \ldots + \binom{h}{h}] = 2^{h-1} s^h p^{(h-1)/2} \bmod q$. Since $h$ is odd and $(p/q) = +1$, this implies $(S/q) = (s/q)$. □

Now let $h$ be the class number of $k = \mathbb{Q}(\sqrt{p}\,)$; since the discriminant of $k$ is a prime, $h$ is odd by genus theory. The factorization above implies $(N + M^2\sqrt{p}\,) = 2\mathfrak{q}^2\mathfrak{a}^4$ for ideals $\mathfrak{q}$ and $\mathfrak{a}$ of norm $q$ and $e$, respectively. Writing $\mathfrak{q}^h = (\lambda)$ and $\mathfrak{a}^h = (\alpha)$ we get

$$(N + M^2\sqrt{p}\,)^h = 2^h \varepsilon_p \lambda^2 \alpha^4$$

with $\varepsilon_p > 1$ as above. Again we form the difference between this equation and its conjugate and then reduce modulo $\overline{\lambda}$. With the help of Lemma 5 we now see that $[\sqrt{p}/\overline{\lambda}] = [\varepsilon_p/\overline{\lambda}]$, and this proves as above that $(q/p)_4 = 1$ is necessary for $\mathcal{T}^{(\psi)}(p)$ to be solvable.

We summarize our discussion in the following

**Theorem 6.** *Consider the elliptic curve $E : y^2 = x(x^2 - 4pq^2)$, where $p \equiv 1 \bmod 8$ and $q \equiv 3 \bmod 4$ are primes such that $(p/q) = +1$. If the torsor $\mathcal{T}^{(\psi)}(p) : N^2 = pM^4 - 4q^2 e^4$ has a rational solution, then $(2/p)_4 = (q/p)_4 = +1$. Moreover, we always have $\mathrm{III}(E/\mathbb{Q})[2] = 0$, and $\mathrm{III}(\widehat{E}/\mathbb{Q})[2]$ has order divisible by $4$ unless possibly when $(2/p)_4 = (q/p)_4 = +1$.*

Our results aren't as complete as those in [10]: if $(2/p)_4 = 1$ and $(q/p)_4 = -1$, for example, then we know that the torsors $\mathcal{T}^{(\psi)}(b_1)$ are not solvable for $b_1 \in \{-1, q, -q\}$, and we also know that one of $\mathcal{T}^{(\psi)}(2)$ and $\mathcal{T}^{(\psi)}(-2)$ is not solvable (assuming the finiteness of $\mathrm{III}(\widehat{E}/\mathbb{Q})$ we can even predict that exactly one of them is solvable), but we can not tell which. It would be interesting to find criteria that allow us to do this.

We still have to provide a proof for

**Lemma 7.** *Let $\pi = a + bi$ be a prime in $\mathbb{Z}[i]$ with norm $p \equiv 1 \bmod 8$, and assume that $4 \mid b$ and $a \equiv 1 \bmod 4$. Then $[\pi/\overline{\pi}]_4 = (-4/p)_8$.*

*Proof.* Since $\pi = a + bi \equiv 2bi \bmod \overline{\pi}$, we find $[\pi/\overline{\pi}]_4 = [2i\overline{\pi}]_4 [b/\overline{\pi}]_4$. Now $[2i\overline{\pi}]_4 = (-4/p)_8$ since $-4 = (2i)^2$, hence it is sufficient to prove that $[b/\overline{\pi}]_4$.

Assume first that $b \equiv 4 \bmod 8$; then $b = 4b'$ for some odd $b'$, and using $[-1/\overline{\pi}]_4 = +1$ we find $[b/\overline{\pi}]_4 = (2/p)[b'/\overline{\pi}]_4 = [b'/\overline{\pi}]_4 = [\overline{\pi}/b']_4 = [a/b']_4$. But quartic residue symbols whose entries are rational integers are trivial (see e.g. [7]), and our claim follows. □

## 4 How did Pépin prove it?

To begin with, it is not clear at all that Pépin really had proofs (in his first paper, he writes "les théorèmes que je propose" which might indicate that he could not prove them); but even if he had I do not believe that he exploited the structure of the class groups of quadratic fields (or quadratic forms), although the necessary techniques had been introduced by Dirichlet [4]. I rather guess that, for proving e.g. Proposition 2, Pépin started with the equation $px^4 - 36y^4 = z^2$ and plugged in the quadratic form for $p$; then one gets $(9a^2 + 4b^2)x^4 - 36y^4 = z^2$, hence

$$9(ax^2 + 2y^2)(ax^2 - 2y^2) = (z + 2bx^2)(z - 2bx^2).$$

Using some form of descent it might be possible to derive a contradiction (Dickson [3, p. 628] mentions a similar method used by Kramer [8]), but I cannot see how. Here is a related proof of Proposition 2 whose idea is taken from Silverman [25, Ch. X]: write $p = 9a^2 + 4b^2$ with $a, b > 0$ and assume that $px^4 - 36y^4 = z^2$ has an integral solution with $x, y, z > 0$. Then it is straightforward to verify the identity

$$p(2bx^2 + 6y^2)^2 = (px^2 + 12by^2)^2 - 9a^2z^2. \tag{1}$$

Now the right hand side factors, and if $x$ is odd, the factors have a factor 2 in common. Since both factors must be positive we get

$$\begin{aligned} px^2 + 12by^2 + 3az &= 2pr^2, & px^2 + 12by^2 - 3az &= 2s^2 \quad \text{or} \\ px^2 + 12by^2 + 3az &= 2r^2, & px^2 + 12by^2 - 3az &= 2ps^2. \end{aligned}$$

If $3 \nmid x$, then both cases are impossible modulo 3. If $3 \mid x$, however, then $3 \mid z$ by the original equation, and this implies $3 \mid r$; since $3 \nmid b$, we conclude that $3 \mid y$, which contradicts the assumption that $(x, y) = 1$. Finally it is easily checked that the case of even $x$ leads to the very same equations.

The identity (1) used in the proof above is easily generalized: for primes $p = \alpha^2 a^2 + 2\beta ab + \gamma b^2$ and $m = \alpha^2\gamma - \beta^2$ as in Theorem 3 one finds

$$mp(\alpha y^2 + bx^2)^2 = (\alpha px^2 + mby^2)^2 - (\alpha^2 a + \beta b)^2 z^2.$$

Yet another possible proof would start by parametrizing the corresponding conic and performing what is nowadays called a second simple 2-descent (see Cremona [2]). In fact, a version of this technique was used already by Euler [5] and thus could have been known to Pépin. Because of its relevance to modern algorithms we will show in the next section that the second descent has one of its roots in the work of Euler.

## 5 The second descent in Euler's papers

One of the earliest papers of Euler on diophantine analysis where he uses Fermat's technique of infinite descent is [5]: in this paper, Euler proves Fermat's

Last Theorem for $n = 4$ (Thm. 1), and he almost shows that the only rational points on $y^2 = x^3 + 1$ are $(-1, 0)$, $(0, \pm 1)$ and $(2, \pm 3)$.

Let us take the time to have a closer look at Euler's proof. He starts by writing $x = a/b$ for coprime integers $a$ and $b$, and multiplying through by $b^4$ he finds that $ba^3 + b^4 = b(a + b)(a^2 - ab + b^2)$ is a square. Putting $a + b = c$, this shows that $bc(c^2 - 3bc + 3b^2)$ is a square. A simple calculation shows that $(b, c) = (b, c^2 - 3bc + 3b^2) = 1$ and $(c, c^2 - 3bc + 3b^2) = (c, 3)$, hence there are two cases to consider. Assume first that $3 \nmid c$; then the three factors $b$, $c$ and $c^2 - 3bc + 3b^2$ are all squares (up to units). Now Euler writes the square root of the last number in the form $\frac{m}{n}b - c$ with coprime integers $m$ and $n$ (this is clearly allowed as long as $b \neq 0$), that is, he puts

$$c^2 - 3bc + 3b^2 = \left(\frac{m}{n}b - c\right)^2. \tag{2}$$

From (2) one easily deduces (using $b \neq 0$) that $\frac{b}{c} = \frac{2mn - 3n^2}{m^2 - 3n^2}$. Again there are two cases $3 \nmid m$ and $3 \mid m$ to consider; in the first case, numerator and denominator in the last equation are coprime, hence $b = 2mn - 3n^2$ and $c = m^2 - 3n^2$ (up to sign, which can be subsumed into $b$ and $c$, however). Thus $b/n^2 = 2\frac{m}{n} - 3$ and $c = m^2 - 3n^2 = (m - \frac{p}{q}n)^2$ for coprime integers $p$ and $q$ (this is Diophantus' substitution again). Collecting these relations shows that $b/n^2 = (p^2 - 3pq + 3q^2)/pq$, and since $b$ is a square, so is $pq(p^2 - 3pq + 3q^2)$. But now Euler is home (the other cases mentioned above can be treated similarly), and a simple reference to descent shows that the only solutions to $y^2 = x^3 + 1$ are those for which $(p, q)$ is not "smaller" than the original solution $(b, c)$. A computation (that Euler skips) shows that this happens only if $b = 1$, which in turn leads directly to the integral points $(-1, 0)$, $(0, \pm 1)$ and $(2, \pm 3)$ on the elliptic curve; Euler only mentions the solution $x = 2$, however.

## Some Comments on Euler's Proof

Next let us discuss Euler's proof from a modern point of view. His substitution $c = a + b$ is easily explained: since $x = \frac{a}{b}$, this amounts to putting $x = X + 1$, that is, to switching from $y^2 = x^3 + 1$ to the model $y^2 = X(X^2 - 3X + 3)$. In other words: Euler chooses a coordinate system in which the rational point of order 2 is in the origin.

Euler's next step is to derive that $b$, $c$ and $c^2 - 3bc + 3b^2$ are all squares if $3 \nmid c$. Writing $c = M^2$ and $b = e^2$, this is equivalent to $N^2 = M^4 - 3M^2e^2 + 3e^4$: but this is nothing but the torsor $\mathcal{T}^{(\psi)}(1)$. In order to show that the only rational point on $\mathcal{T}^{(\psi)}(1)$ satisfies $e = 0$, Euler considers the underlying conic: in fact, his substitution (2) is a trick he learned from Diophantus: what Euler is doing here (apparently without being aware of the geometrical background) is parametrizing the conic $y^2 = c^2 - 3bc + 3b^2$ by using the line $y = \frac{m}{n}b - c$ with the rational slope $\frac{m}{n}$ through the rational point $b = 0$, $y = -c$.

## 6 The case $n = 7$ of Fermat's Last Theorem

Although it is known that there exists an elementary proof of Fermat's Last Theorem for $n = 7$, most people seem to think that this proof is utterly complicated. In Ribenboim's book [23, pp. 46–67], one finds Lamè's proof from 1839 as simplified by Lebesgue up to the point where one is left with showing that a certain quartic has only trivial solutions. This quartic appears in all the elementary proofs and is the subject of various simplifications (e.g. by Pépin [16] and Genocchi [6]; see Nagell [13]). Usually one uses a descent à la Fermat to disprove the existence of nontrivial rational points, but since the quartic in question is an elliptic curve (it has conductor $7^2$) with a rational point of order 2, we can use the machinery presented in Section 3.

Suppose that $(x, y, z)$ is a rational solution of the Fermat equation for $n = 7$. Form the polynomial $X^3 - pX^2 + qX - r = (X-x)(X-y)(X-z)$; then $p, q$ and $r$ are rational numbers (even integers), and Newton's formulas (see e.g. [12]) give

$$x^7 + y^7 + z^7 = p^7 - 7p^5q + 7p^4 + 14p^3q^2 - 21p^2qr - 7pq^3 + 7pr^2 + 7q^2r.$$

Replacing $r$ by $pq - r$ gives $x^7 + y^7 + z^7 = p^7 - 7p^4r + 7p^2qr + 7pr^2 - 7q^2r$. Since $x^7 + y^7 + z^7 = 0$, the right hand side must vanish:

$$p^7 - 7p^4r + 7p^2qr + 7pr^2 - 7q^2r = 0 \qquad (3)$$

Suppose first that $p \neq 0$. Then we can substitute $q = p^2Q$ and $r = p^3R$, and after canceling $p^7$ we find $7R^2 - 7R(1 - Q + Q^2) + 1 = 0$. Since this equation has rational solutions, its discriminant must be a square, and we get $(Q^2 - Q + 1)^2 - 4/7 = \square$. Now we write $2Q - 1 = \frac{s}{t}$ and find

$$u^2 = s^4 + 6s^2t^2 - \frac{1}{7}t^4. \qquad (4)$$

At this point, the published proofs start performing some descent on (4) to show that it has only the obvious solutions. But modern eyes recognize (4) at once as an elliptic curve! Divide by $t^4$ and put $U = u/t^2$ and $S = u/t$. Then $(U - S^2 - 3)(U + S^2 + 3) = -\frac{64}{7}$, and if you write $T = U - S^2 - 3$, then $U + S^2 + 3 = -\frac{64}{7T}$. Substracting the last two equations gives you $T + \frac{64}{7T} = -2S^2 - 6$. Multiply by $T^2$ and put $Y = ST$; then $T^3 + 6T^2 + \frac{64}{7}T = -2Y^2$. Now we're almost there: multiplying by $-2^{-3} \cdot 7^6$ and putting $y = \frac{1}{2}7^3Y$ and $x = -\frac{1}{2}7^2T$ we find

$$E : y^2 = x(x^2 - 3 \cdot 7^2x + 2^4 \cdot 7^3)$$

(you will find that we just applied the algorithm in the book of Cassels [1]). The theorem of Nagell-Lutz shows that $E(\mathbb{Q})_{\text{tors}} = \{\mathcal{O}, (0, 0)\}$, and a simple 2-descent shows that $E$ has rank 0. Transforming back we find that $t = 0$ is the only solution of the quartic (4), which in turn means that (3) has no rational solution with $p \neq 0$.

If $0 = p = x + y + z$, on the other hand, then we get $x^7 + y^7 = (x+y)^7$, hence
$$0 = x^7 + y^7 - (x+y)^7 = 7xy(x+y)(x^2 + xy + y^2)^2.$$

Thus any solution of the Fermat equation for $n = 7$ with $xyz \neq 0$ must satisfy $x^2 + xy + y^2 = 0$ which is impossible in integers.

# References

[1] J.W.S. Cassels, *Lectures on elliptic curves*, Cambridge University Press, 1991

[2] J. Cremona, *Higher descents on elliptic curves*, preprint 1998

[3] Dickson, *History of the Theory of Numbers II*, Chelsea Publishing, New York 1952

[4] L. Dirichlet, *Untersuchungen über die Theorie der quadratischen Formen*; Abh. Königl. Preuss. Akad. Wiss. 1833, 101–121; Werke I, 195–218

[5] L. Euler, *Theorematum quorundam arithmeticorum demonstrationes*, Comm. Acad. Sci. Petrop. **10** (1738) 1747, 125–146; Opera Omnia Ser. I vol. II, Commentationes Arithmeticae, 38–58

[6] A. Genocchi, *Géneralization du théorème de Lamé sur l'impossibilité de l'équation $x^7 + y^7 + z^7 = 0$*, C. R. Acad. Sci. Paris **82** (1876), 910–913

[7] K. Ireland, M. Rosen, *A classical introduction to modern number theory*, Springer Verlag, 2nd. ed. 1990

[8] A.E. Kramer, *De quibusdam aequationibus indeter. quarti gradus*, Diss. Berlin 1839

[9] F. Lemmermeyer, *Rational quartic reciprocity*, Acta Arith. **67** (1994), 387–390

[10] F. Lemmermeyer, *On Tate-Shafarevich groups of some elliptic curves*, Proc. Conf. Graz 1998

[11] C.-E. Lind, *Untersuchungen über die rationalen Punkte der ebenen kubischen Kurven vom Geschlecht Eins*, Diss. Univ. Uppsala 1940

[12] P. Morandi, *Field and Galois Theory*, Springer-Verlag 1996

[13] T. Nagell, *Introduction to number theory*, John Wiley & Sons (1951)

[14] E. Netto, *Review 06011301*, Jahrbuch der Fortschritte der Mathematik **6** (1874), p. 113

[15] Th. Pépin, *Théorèmes d'analyse indéterminée*, C. R. Acad. Sci. Paris **78** (1874), 144–148

[16] Th. Pépin, *Impossibilité de l'équation $x^7 + y^7 + z^7 = 0$*, C. R. Acad. Sci. Paris **82** (1876), 676–679, 743–747

[17] Th. Pépin, *Théorèmes d'analyse indéterminée*, C. R. Acad. Sci. Paris **88** (1879), 1255–1257

[18] Th. Pépin, *Nouveaux théorèmes sur l'équation indéterminée $ax^4 + by^4 = z^2$*, C. R. Acad. Sci. Paris **91** (1880), 100–101

[19] Th. Pépin, *Nouveaux théorèmes sur l'équation indéterminée $ax^4 + by^4 = z^2$*, C. R. Acad. Sci. Paris **94** (1882), 122–124

[20] L. Rédei, *Die Diophantische Gleichung $mx^2 + ny^2 = z^4$*, Monatsh. Math. Phys. **48** (1939), 43–60

[21] L. Rédei, H. Reichardt, *Die Anzahl der durch 4 teilbaren Invarianten der Klassengruppe eines beliebigen quadratischen Zahlkörpers*, J. Reine Angew. Math. **170** (1933), 69–74

[22] H. Reichardt, *Einige im Kleinen überall lösbare, im Großen unlösbare diophantische Gleichungen*, J. Reine Angew. Math. **184** (1942), 12–18

[23] P. Ribenboim, *13 lectures on Fermat's Last Theorem*, Springer-Verlag (1979)

[24] H.E. Rose, *On some classes of elliptic curves with rank two or three*, Univ. Bristol Math. Res. Report PM–97–01

[25] J. Silverman, *The arithmetic of elliptic curves*, Graduate Texts in Mathematics, Springer-Verlag, 1986